\documentclass[12pt]{amsart}
\usepackage{amsmath,amssymb,amsfonts,amsthm,amsopn}

\newcommand{\tfa}{time-frequency analysis}

\newcommand{\stft}{short-time Fourier transform}

\newcommand{\tf}{time-frequency}

\newcommand{\tfs}{time-frequency shift}

\newcommand{\aaf}{A_a^{\varphi_1,\varphi_2}}

\newcommand{\modsp}{modulation space}
\newcommand{\psdo}{pseudodifferential operator}

\newtheorem{tm}{Theorem}[section]
\newtheorem{lemma}[tm]{Lemma}

\newcommand{\rem}{\textsl{REMARK:}}

\newtheorem{theorem}{Theorem}[section]
\newtheorem{corollary}[theorem]{Corollary}
\newtheorem{definition}[theorem]{Definition}

\newtheorem{proposition}[theorem]{Proposition}

\theoremstyle{remark}

\newcommand{\beqa}{\begin{eqnarray*}}
\newcommand{\eeqa}{\end{eqnarray*}}

\DeclareMathOperator*{\supp}{supp}

\newcommand{\field}[1]{\mathbb{#1}}
\newcommand{\bR}{\field{R}}        
\newcommand{\bN}{\field{N}}        
\newcommand{\bZ}{\field{Z}}        
\def\la{\lambda}

\def\om{\omega}

 \def\cS{\mathcal{S}}
 \def\cD{\mathcal{D}}

 \def\cM{\mathcal{M}}

 \def\cC{\mathcal{C}}
 
 \def\cO{\mathcal{O}}

\def\a{\aleph}

\def\vgf{V_gf}

\def\rd{\bR^d}

\def\rdd{{\bR^{2d}}}

\def\lrd{L^2(\rd)}
\def\lrdd{L^2(\rdd)}

\def\mvv{M_v^1}

\def\intrdd{\int_{\rdd}}

\def\R{\right)}

\def\<{\left<}
\def\>{\right>}

\def\inv{^{-1}}

\def\mv1{M_v^1}

\def\Lmpq{L_m^{p,q}}

\def\mpq{M^{p,q}}
\def\Mmpq{M_m^{p,q}}
\def\phas{(x,\omega )}

\def\c{\hfill\break}
\def\ni{\noindent}

\def\o{\omega}
\def\a{\alpha}
\def\b{\beta}

\def\R{\mathbb{R}}
\def\Ren{\mathbb{R}^d}
\def\Renn{\mathbb{R}^{2d}}

\def\Qs{{Q_s}}

\def\sch{\mathcal{S}}

\def\Fur{\mathcal{F}}

\def\f{\varphi}

\def\gaw{A_a^{\f_1,\f_2}}

\def\Sn2{S_{2}(L^{2}(\Ren))}
\def\S1{S_{1}(L^{2}(\Ren))}
\def\sig00{\sigma_{0,0}}

\def\la{\langle}
\def\ra{\rangle}

\def\ts{\tau_s}

\def\wei{\la\cdot\ra}

\newcommand{\vf}{\varphi}
\newcommand{\win}{\vf _0}
\newcommand{\speci}{\cS _{\cC }}
\begin{document}
\begin{abstract}
We study the composition of   \tf\ localization operators
(wavepacket operators)   and develop a  symbolic calculus of such
operators on  modulation spaces. The use of  \tf\ methods
(phase space methods) allows the   use of  rough symbols of
ultra-rapid growth in place of smooth symbols in the standard classes. 
As the main application it is shown that, in general, a localization
operators possesses the Fredholm property, and thus  its  range is closed
in the target space. 
\end{abstract}

\title[Symbolic calculus for localization operators]{Symbolic calculus and
Fredholm property \\ for localization operators}
\author{Elena Cordero and Karlheinz Gr\"ochenig}
\address{Department of Mathematics,  University of Torino, Italy}

\address{Faculty of Mathematics, University of Vienna, Nordbergstra\ss e 15,
A-1090 Vienna, Austria}
\email{cordero@dm.unito.it,karlheinz.groechenig@univie.ac.at}
\maketitle

\section{Introduction}

By a symbolic calculus is meant  a mapping from a parameter space of symbols to
 a class of operators and the investigation of this functional
 dependence. The prototype of a symbolic calculus  is the symbolic
 calculus of \psdo s. This classical  symbolic calculus has lead to an
 understanding of the composition of \psdo s, the construction on
 an approximate inverse, a  so-called parametrix, and the regularity
 properties of partial differential operators~\cite{Hormander3,Shubin91,taylor81}.

In this paper we study a different kind of  symbolic calculus, namely
for (time-frequency) localization operators. This class of operators
occurs in various 
branches of mathematics under such names as Toeplitz operators, 
(anti) Wick operators,  time-frequency multipliers, and others. Their applications range from
quantization procedures (Berezin quantization) in quantum mechanics, via signal analysis, to the 
approximation of \psdo s. 
Localization operators form a special case of the  wave-packet
operators of  Cordoba and
Fefferman~\cite{CF78} and  have been used to approximate general \psdo
s. Technically,  localization operators  may be regarded as  a special   
class of \psdo s and may  therefore be investigated with 
pseudodifferential calculus. This approach has led to many interesting
results in \cite{PaoloToft,PT04,Toft04,Toftweight}.

 Our approach to localization operators is based exclusively on \tf\
 methods (phase-space methods). In view of the very definition of this
 class of operators  below, this approach is   not only
 natural, but  it also leads to  very strong results with respect to
 the admissible symbols and ``window functions''. In  treating the
 composition of two localization operators one of the symbols may be
 rough (in local $L^\infty $) and even possess ultra-rapid growth. 

Whereas the product of 
two \psdo s is again a \psdo\ (there are many algebras of \psdo s),
the composition of two localization 
operators is not a localization operator in general. 
This additional difficulty has  captured the interest  of several
authors and has generated some remarkable ideas. 
 An exact product formula for localization operators is
presented in \cite{du-wong-product}. However, since it works only under
very restrictive conditions  and is unstable, it is not amenable to a
detailed  analysis of mapping or compactness properties. Therefore
many authors resort to asymptotic expansions that realize the
composition of two localization operators as a sum of localization
operators and a controllable
remainder~\cite{Ando-Morimoto02,CR03,Lerner01,Lernercubo,Tataru02}.
These contributions  were mainly  motivated by PDEs and energy
estimates, and therefore use smooth symbols  that are defined by
differentiability properties, such as the traditional H\"ormander or
Shubin classes, and Gaussian windows. For applications in quantum mechanics and signal
analysis, alternative notions of smoothness --- ``smoothness in
phase-space'' or quantitative measures of ``\tf\ concentration'' ---
have turned out to be useful. This point of view is pursued  in
several  recent 
investigations of localization operators and usually involves  modulation
spaces (see~\cite{CR03,Toft04,Toftweight}  and  references therein).
 In the context of \modsp s, much rougher symbols and more general 
``window functions'' can be used for localization operators than have
been  considered in the  studies~\cite{CG02,Daube88,RT93}.

Modulation spaces --- though still not as well known as standard
smoothness spaces --- are the appropriate function spaces for \tfa ,
and in several cases have been shown to furnish optimal results. For
instance, they arise as the optimal symbol classes in the study of
boundedness and Schatten class properties of localization operators~\cite{CGAMS04,CG02}

Our analysis of the symbolic calculus for localization operators is
based on a composition formula in \cite{Ando-Morimoto02} and (in full
generality) in~\cite{CR03}. It  was formulated for
symbols in certain Shubin classes and windows in the Schwartz class. 
This calculus  seems  most suitable for understanding  the composition of
localization operators with respect to their action on phase-space
(the \tf\ plane).  

\vspace{3 mm}

\textbf{The Short-Time Fourier Transform and Localization Operators.}
To be more specific and to formulate our results, we first define the \stft\ and  
introduce the class of  localization operators. 

 The operators of translation and modulation are defined to be 
\begin{equation}
  \label{eqi1}
 T_xf(t)=f(t-x)\quad{\rm and}\quad M_{\o}f(t)= e^{2\pi i \o
  t}f(t) \, .
\end{equation}
We often combine translations and modulations into \emph{\tf\, shifts}
(phase-space shifts in physical terminology).  Write  $z=\phas \in
\rdd$, then  the general \tfs\  is
defined by 
\begin{equation}\label{pi}
\pi(z)=M_\om T_x.
\end{equation}

Associated to \tfs\ is an important transform, the \emph{short-time
  Fourier transform} (STFT), which is also called coherent state
transform, Gabor 
transform, windowed Fourier transform, and the like. The STFT of a
function or distribution $f$ with respect to a fixed non-zero
\emph{window} function  $g$ is given by
 \begin{equation}
   \label{eqi2}
   V_gf(x,\o)=\int_{\Ren}
 f(t)\, {\overline {g(t-x)}} \, e^{-2\pi i\o t}\,dt = \la f,M_\o T_x
 g\ra = \langle f, \pi (z) g\rangle \, , 
 \end{equation}
whenever the integral or the bracket $\langle  \cdot , \cdot \rangle
$ is well-defined. 
 The \stft\   can be  defined on
many pairs of  distribution spaces and test functions.  For instance,
$V_gf$, as a sesquilinear 
form,  maps $L^2(\rd
) \times L^2(\rd )$ into $\lrdd $ and $\sch(\Ren)\times\sch(\Ren)$
into $\sch(\Renn)$. Furthermore, $V_gf$  can be extended  to a map
from $\sch'(\Ren)\times\sch'(\Ren)$ into
$\sch'(\Renn)$~\cite[p.~41]{book}. The \stft\ is  the appropriate
tool   for defining localization operators.

Next, let $a$ be a  \emph{symbol} on the \tf\
``plane'' $\Renn$ and  choose two windows $\f _1, \f _2$ on $\rd$,
then the  localization operator $\aaf $  is defined   as
\begin{equation}
  \label{eqi4}
\aaf f(t)=\int_{\Renn}a \phas V_{\f _1}f \phas M_\omega T_x \f _2 (t)
\,
dx d\omega \, .
\end{equation}
Taking the ``inner product'' with a test function $g$, the definition
of  $\aaf$ can be written in a weak sense,  namely,
\begin{equation}\label{anti-Wickg}
 \la \gaw f,g\ra=\la aV_{\f_1}f, V_{\f_2}g\ra=\la
 a,\overline{V_{\f_1}f}\,  V_{\f_2}g\ra.
\end{equation}

If $a \in \cS '(\Renn )$ and $\f _1, \f _2 \in \cS (\rd )$, then
\eqref{anti-Wickg} defines a  continuous operator from $\cS (\rd
)$ to $\cS ' (\rd )$. If $\f _1(t)  = \f _2 (t) = 2^{d/4}e^{-\pi
t^2}$, then $A_a = \aaf $ is well-known as  (anti-)Wick operator
and the mapping $a \to \aaf  $ is interpreted as a quantization
rule \cite{Berezin71,du-wong-product,Lerner01,Shubin91,Wong02}.
Furthermore, this definition is also a special case of the wave packet
operators of Cordoba and Fefferman~\cite{CF78,folland89}, and $\gaw$  serves  as an
approximation of the \psdo\ $a(x,D)$. 

While previous work -- to a large extent -- uses  a combination of
methods and is always focused on the model of \psdo s, we will use
exclusively \tfa\ in our study of localization operators. In a sense,
this is perfectly natural, because the definition of a localization
operators is in terms of the basic transform of \tfa , namely the
STFT.

\vspace{3 mm}

\textbf{Results.}
Our  starting
point is the following composition formula for two  localization
operators derived in~\cite{CR03}: 
\begin{equation}\label{c1}
\gaw A^{\f_3,\f_4}_b= \sum_{|\a|=0}^{N-1} \frac{(-1)^{|\a|}}{\a!}
A^{\Phi_\a,\f_2}_{a\partial^{\a} b}+E_N
\end{equation}

The essence of this formula is that the product of two localization
operators can be written as a sum of localization operators and a
remainder term, which  is ``small''. 

In the spirit of the  classical symbolic calculus, this formula was
derived in \cite[Thm. 1.1]{CR03} for  \emph{smooth} symbols belonging
to some Shubin class $S^m (\rdd )$ and for
windows in the  Schwartz class $\cS (\rd )$. In this case the remainder term
is regularizing and maps $\lrd $ into $\cS (\rd )$.

The goal of this paper is much more ambitious. We will establish the
validity of~\eqref{c1} on the \modsp s (Theorem~\ref{mainA} and
Proposition~\ref{mainadj}). These  function spaces are
defined by imposing a weighted mixed $L^p$-norm on the STFT, e.g.,
$\|f\|_{M^p} = \|V_g f \|_{p}$ (for suitable  window $g$). In view of
definition~\eqref{anti-Wickg} \modsp s  are bound  to arise
in the   in-depth treatment of the symbolic calculus.

The extension of the composition formula ~\eqref{c1}  is far from a
routine generalization  and requires the entire arsenal of \tfa . 
The derivation of a symbolic calculus  is necessarily technical, and 
 therefore we would like to   point out explicitly  the innovative
 features and the level of generality of our results. 

(i) \emph{Rough symbols.} While in \eqref{c1} the symbol $b$ must be
$N$-times differentiable, the symbol $a$ only needs to be locally
bounded. The classical results in symbolic calculus require both
symbols to be smooth.

(ii) \emph{Growth conditions on symbols.} The symbolic calculus in
\eqref{c1} can handle symbols with ultra-rapid growth (as long as it
is compensated by a decay of $b$ or vice versa). For instance, $a$ may grow
subexponentially as $a(z) \sim e^{\alpha |z|^\beta }$ for $\alpha >0$
and $0< \beta <1$. This goes far beyond the usual polynomial growth and
decay  conditions. 

(iii) \emph{General window classes.} We provide a precise description
of the admissible windows $\f _j$ in \eqref{c1}. Usually only the
Gaussian $e^{-\pi x^2}$ or Schwartz functions are considered as 
windows~\cite{Ando-Morimoto02,CR03,Lerner01}.  

(iv) \emph{Size of the remainder term.} We derive norm estimates for
the size of the remainder term $E_N$ that depend explicitly on the
symbols $a,b$ and the windows $\f _j$. In applications to PDE, e.g.~\cite{Hormander3,taylor81}
it is important that the remainder be regularizing, but this property
does not exclude the possibility that $E_N$ is large in
norm. Therefore 
norm estimates  provide important additional information, see
also~\cite{Lerner01,Lernercubo} for more motivation.

\vspace{3 mm}

\textbf{The Fredholm Property of Localization Operators.} By choosing
$N=1$, $\f _1 = \f _2 = \f$ with 
$\|\f \|_2=1$, $a (z) \neq 0$, and $b = 1/a$, the composition
formula~\eqref{c1} yields the following important special case:
\begin{equation}
  \label{eq:mu1}
  A^{\f, \f } _a  A^{\f, \f } _{1/a} = A_1^{\f , \f } + R = \mathrm{I}
  + R \, .
\end{equation}
Under mild conditions on $a$ we will show that $R$ is  compact  and
that $ A^{\f, \f } _a $ is a Fredholm operator between two \modsp s $\mpq $
and $\Mmpq $ (with different weights), see Theorem~\ref{F1}.   This 
theorem is remarkable because it  works even for ultra-rapidly
growing symbols such as $a(z) = e^{\alpha |z|^\beta } $ for $\alpha
>0$ and $0<\beta <1$. For comparison, the reduction of localization
operators to standard pseudodifferential calculus requires 
  elliptic or hypo-elliptic symbols,   and the proof of the Fredholm
property works only under severe restrictions, see~\cite{PaoloToft}. 

\vspace{3 mm}

\textbf{Itinerary of the Proof of the Symbolic Calculus.} To give the
reader some insight how  the expansion formula~\eqref{c1} is derived,
we sketch the main  arguments developed in~\cite{CR03}.

The composition of two localization operators can be formally written
as 
\begin{equation}\label{start}
  \gaw A^{\f_3,\f_4}_{b}f =\int_{\Renn}\int_{\Renn}a (y)b(z)
  V_{\f_3}f(z)\langle \pi(z)\f_4 , \pi (y) \f _1\rangle \,
  \pi(y)\f_2 \,dydz\, .  
\end{equation}
Assuming the symbol $b$ to be $N$-times differentiable, we  expand
$b(z)$ into a Taylor series around $y$ and obtain 
$$b(z)=\sum_{|\a|\leq N} \partial^{\a} b
(y)\frac{(z-y)^{\a}}{\a!}+b_N(y,z)\, ,$$
 where the remainder is given by  
$$b_N(y,z)=N\sum_{|\a|=N}\int_0^1(1-t )^{N-1}\partial^\a
b(y+t (z-y))\,dt \,\,  \frac{(z-y)^\a}{\a!}\, .
$$
For each $\alpha $ with $|\alpha | < N$ we obtain a term of the form 
$$
\frac{1}{\alpha ! }\, \int _{\Renn } a(y) \, \partial ^\alpha b(y) \, \Big(
\int _{\Renn } ( z-y )  ^\alpha V_{\f _3 } f(z)  \, 
\langle \pi (z) \f _4 , \pi (y) \f _1 \rangle \, dz \Big) \pi
(y) \f _2 \, dy \, .
$$
The inner integral over $z$ is linear in $f$ and ``covariant'' in $y$,
and can be  expressed as an STFT $\langle f, \pi (y) \Phi
_\alpha \rangle $, where $\Phi _\alpha $ depends on $\f _1, \f _3, \f
_4$. Using formulas for the moments of the STFT, this new window was 
calculated explicitly in~\cite{CR03}  to be 
\begin{equation}\label{newwindc}
\Phi_\a=\frac1{(2\pi i)^{|\a_2|}}\sum_{\b\leq\a}
\binom{\a}{\b}(-1)^{|\b_1|}\la \f_3,
X^{\a_1-\b_1}\partial^{\a_2-\b_2}\f_4\ra
X^{\b_1}\partial^{\b_2}\f_1 \, .
\end{equation}
[Note that if $\f _1, \f _3, \f_4 \in \cS (\Ren )$, then also $\Phi
_\alpha \in \cS (\Ren )$.] Consequently, for the terms $|\alpha | <
N$, we obtain the  localization operators
$$ 
\frac{1}{\alpha !}\int _{\Renn }   a(y) \partial ^\alpha b(y)  \,
V_{\Phi _\alpha } f (y)  \, \pi (y) \f _2 \, dy = \frac{1}{\alpha !}
A^{\Phi _\alpha , \f 
  _2} _{a\partial ^\alpha b} f \, .
$$
The terms corresponding to $|\alpha | = N$ can be collected to a
remainder $E_N$ given informally by 
\begin{equation}
  \label{eq:neu10c}
  E_N f = \intrdd \intrdd a(y) b_N (y,z) 
  \, V_{\f _3} f (z)  \langle   \pi (z) \f _4, \pi (y) \f _1\rangle
  \pi (y) \f _2 \, dy dz  \, . 
\end{equation}
By summing over $\alpha , |\alpha | \leq N$, we obtain the expansion
formula~\eqref{c1}. In~\cite[Thm. 1.1]{CR03} this formal idea was made
rigorous for symbols $a,b$ contained  in a Shubin class and 
windows  $\f _j $ in the Schwartz class.

For the  extension of the symbolic calculus~\eqref{c1} to rough symbols and general
modulation spaces we need to check  that all terms in \eqref{c1} are
well-defined on \modsp s and then  derive norm estimates for their
size.



To carry out this plan, we will proceed along the following steps. In Section~2 
 we introduce the modulation spaces and list some of their main
 properties. In Section~3  we
present further  \tf \, tools, such as properties of the STFT.  In
Section~4  we investigate the mapping properties of localization
operators between \modsp s. These results enable us to give a rigorous
meaning to most  terms in the symbolic calculus~\eqref{c1}. Section
$5$ is devoted to a careful analysis  
of the remainder term $E_N$,  and we will formulate conditions for the
boundedness and compactness of  $E_N$ on \modsp s.  This  is
the most technical part of the paper. In  Section~$6$ we  combine the
entire machinery and prove the symbolic calculus~\eqref{c1}  in
  the framework of \modsp s.  In the  last Section $7$ we study  the
  Fredholm property of  localization operators. This property  is the main
application of the symbolic calculus  and --- at least for us ---
justifies   its many subtle   technicalities. 

\vspace{3 mm}

\textbf{Notation.} We define $t^2=t\cdot t$, for $t\in\Ren$, and
$xy=x\cdot y$ is the scalar product on $\Ren$. Given a vector
$x=(x_1,\dots,x_d)\in\rd$, the partial derivative with respect to
$x_j$ is denoted by  $\partial _j = \frac{\partial}{\partial x_j}$.
  Given a multi-index
$\a=(\a_1,\dots,\a_d)\geq 0$, i.e., $\alpha \in\bZ^d$ and $\alpha
_j \geq 0$, we write
$\partial^\a=\partial^{\a_1}_1\cdots\partial^{\a_d}_d$; moreover,
we denote by $X^\alpha$ the operator of multiplication: $(X^\alpha
f)(t)=(t^{\a_1}_1\cdots t^{\a_d}_d)f(t)$. 

Given a set $K\subseteq \rdd $, we denote by $meas(K)$ the Lebesgue measure of $K$ and
by  $\chi_K$  the
characteristic  function of  $K$.
  We use the brackets
$\la f,g\ra$ to denote the extension to $\sch ' (\Ren)\times\sch
(\Ren)$ of the inner product $\la f,g\ra=\int f(t){\overline
{g(t)}}dt$ on $L^2(\Ren)$. The space of smooth functions with
compact support on $\rd$  is  denoted by $\cD(\rd)$. 
We denote by $L^0(\rdd ) $ the space  of all  $f\in L^\infty (\rdd ) $
that vanish at infinity, i.e.,  for all $\varepsilon >0$
there exists  $R=  R(\varepsilon) >0$ such that 
$$|f(z)|<\varepsilon,\quad{\rm \,\,for\, a.~a.~}\,|z|>R\, .$$
The Fourier transform is
normalized to be ${\hat
  {f}}(\o)=\Fur f(\o)=\int
f(t)e^{-2\pi i t\o}dt$.

Throughout the paper, we shall use the notation $A\lesssim B$ to
indicate $A\leq c B$ for a suitable constant $c>0$, whereas $A
\asymp B$ means that $c\inv A \leq B \leq c A$ for some $c\geq 1$. The
symbol $B_1 \hookrightarrow B_2$ denotes the continuous embedding of
the linear space $B_1$ into $B_2$.

\section{ Modulation   Spaces }\label{fspaces}

\subsection{Weight Functions}

 For the quantitative description of  decay
properties, we use  weight
functions  on the \tf\ plane. In the sequel $v$ will always be a
continuous, positive,  even, submultiplicative  weight function
(in short, a submultiplicative weight), i.e., $v(0)=1$, $v(z) =
v(-z)$, and $ v(z_1+z_2)\leq v(z_1)v(z_2)$, for all $z,
z_1,z_2\in\Renn.$
We furthermore impose the GRS-condition
\begin{equation}
  \label{eq:2}
  \lim _{n\to \infty } v(nz) ^{1/n} = 1, \quad \quad \forall z \in \rdd
  \, .
\end{equation}
The GRS-condition (introduced by
Gelfand-Raikov-Shilov~\cite{gelfandraikov}) quantifies  the subexponential
growth of $v$ in a precise
manner.
Every weight of the form $v(z) =   e^{a|z|^b} (1+|z|)^s \log ^r(e+|z|)
$ for parameters $a,r,s\geq 0$, $0\leq b < 1$ satisfies the
GRS-condition, whereas the exponential weight $v(z) = e^{a|z|}, a>0,$ does
not. Finally, we assume that $v$ satisfies  the  property
\begin{equation}\label{*}
\int_{0}^{1} v(tz)\,dt \lesssim v(z).
\end{equation}

Associated to every submultiplicative weight we consider the class of
so-called  {\it
  v-moderate} weights $\cM _v$. A  positive, even
weight function  $m$ on $\Renn$ belongs to $\cM _v$ if it  satisfies
the condition
$$
 m(z_1+z_2)\leq Cv(z_1)m(z_2)  \quad   \forall z_1,z_2\in\Renn \, .
$$
 We note that this definition implies that
$\frac{1}{v} \lesssim m \lesssim v $,  $m \neq 0$ everywhere, and that
$1/m \in \cM _v$.

 We will often   use the polynomial weights $\langle \cdot \rangle ^s$    defined by
\begin{equation}
 \la \phas \ra^s= \la z\ra^s=(1+x^2+\o^2)^{s/2},\quad
   z=(x,\o)\in\Renn \,\quad s\in\bR\,  \label{eqc1}
 \end{equation}
and the product weights $\langle \cdot \rangle ^s \otimes \langle
\cdot \rangle ^s$ signifying $\langle x\rangle ^s \langle \omega
\rangle ^s$.  

\subsection{Modulation Spaces}
\label{modspdef}

\hfill\break
\ni
The modulation space norms are a  measure of
the joint time-frequency distribution of a function or distribution.
They occur in many problems in \tfa\ and play an increasing role as
symbol classes for various types of \psdo s.
 For their
basic properties we refer, for instance, to \cite[Ch.~11-13]{book} and
the original literature quoted there.

Let $\win (t)  = e^{-\pi t^2}$ be the Gaussian and let
\begin{equation}
  \label{eq:cn1a}
  \cS _{\cC} = \mathrm{span}\, \{ \pi (z) \win : z\in \rdd \}
\end{equation}
be the linear space of all finite linear combination of \tfs s of the
Gaussian (the ``space of special windows'').

\begin{definition}
 For any moderate weight $m\in \cM _v$ and $1\leq p,q\leq
\infty$ we define the \modsp\ norm on $\cS _{\cC }$ by
$$
\|f\|_{M^{p,q}_m}=\|V_{\vf _0}f\|_{L^{p,q}_m}=\left(\int_{\Ren}
  \left(\int_{\Ren}|V_{\vf _0} f(x,\o)|^pm(x,\o)^p\,
    dx\right)^{q/p}d\o\right)^{1/q}  \, .
$$
If $p,q< \infty $, the \modsp\ $\Mmpq $ is the norm completion of
$\speci $ in the $\Mmpq $-norm. If $p=\infty $ or $q=\infty$, then
$\Mmpq $ is the completion of $\speci $ in the $\sigma (\speci,
\speci )$-topology. Then by definition,   $\Mmpq (\Ren )$ is a Banach space.
 If $p=q$, we write $M^p_m$ instead of
$M^{p,p}_m$, and if $m(z)\equiv 1$ on $\Renn$, then we write
$M^{p,q}$ and $M^p$ for $M^{p,q}_m$ and $M^{p,p}_m$.
\end{definition}

This definition can be read on several levels of generality.
If $v$ and therefore $m\in \cM _v $ grow polynomially, $m(z) = \cO
(|z|^N)$ for some $N\geq 0$, then $\Mmpq $  consists of all tempered
distributions $f\in\sch'(\Ren)$ such that $V_{\win }f\in L^{p,q}_m(\Renn )$
with the  norm $\|f\|_{M^{p,q}_m}=\|V_{\vf _0} f\|_{L^{p,q}_m}$.
If $p,q\leq 2$ and $m\geq 1$, then $\Mmpq $ is  a subspace of
$\lrd $.

However, if $v$ and $m$ grow faster than polynomially, we may need to
appeal to the theory of ultra test functions and ultradistributions to
get a concrete definition of $\Mmpq $ in place of the equivalent
definition as an abstract norm completion~\cite{bjork66,fg89jfa,PT04}. For instance, if $v(z) =
e^{a|z|^\gamma }$ for $\gamma <1$, then $\mvv \subseteq \cS (\rd ) \subseteq
 \cS ' (\rd ) \subseteq M^\infty _{1/v}$, and the elements in
 $\Mmpq $ are those ``ultra distributions'' whose STFT is in
 $L^{p,q}_m$.

We list the some of the main properties of \modsp
s~\cite[Ch.~11]{book}. 

\begin{tm}
  \label{modprop1}
Let $1\leq p,q \leq \infty $,  $v$ be a submultiplicative weight as
above,  and $m \in \cM _v$.

(i)   Then $\Mmpq (\Ren )$ is a Banach space, and
\begin{equation}
  \label{eq:dec1}
  \|\pi (z) f \|_{\Mmpq } \lesssim v(z) \, \|f\|_{\Mmpq } \, .
\end{equation}

(ii) Duality: If $1\leq p,q < \infty $ and $p'= \frac{p}{p-1} $ is
 the conjugate exponent, then
 $(M^{p,q}_m)^*=M^{p',q'}_{1/m}$.
\end{tm}

\rem\ The \modsp s $M^{\infty ,1} _m $ is the dual of $M^{1,0}_{1/m}$
where $M^{1,0}_{1/m}$ is the closure of $\cS _{\cC } $ in the
$M^{1,\infty }_{1/m}$-norm~\cite{BGHO04}. Similarly,  $M^{1, \infty
} _m  =( M^{0,1}_{1/m} )^*$. Thus duality and weak-$^*$ arguments can be used
for all indices $p,q$.

The following statement on norm equivalence is crucial and  will be used
repeatedly. Again $\win(t)=e^{-\pi t^2}$.

\begin{tm}\label{convV*3}
 (i) \cite[Lemma~12.1.1]{book} If  $f,g \in \mvv $, then  
\begin{equation}
  \label{eq:cn10}
  \|V_g f \|_{L^1_v } \, \lesssim \, \|f \|_{M^1_v } \, \|g\|_{M^1_v } \, .
\end{equation}
 (ii)  Equivalent norms~\cite[Thm.~11.4.2]{book}: If $g\in \mvv $ and $f\in \Mmpq$, then
\begin{equation}
  \label{eq:ch1}
 \|V_{g}f \|_{\Lmpq } \lesssim \|g \|_{M^1_v } \, \|f\|_{\Mmpq }
 \lesssim \|V_g f \|_{\Lmpq } \, .
\end{equation}
\end{tm}

The estimates in~\eqref{eq:ch1} imply  that  the definition
of $\Mmpq $ is  independent of  the window 
$g$ and that  $\|\vgf \|_{L^{p,q}_m}$ is an equivalent norm for any
$g\in M^1_v$. To understand the various conditions on the occurring
windows, it is helpful to  keep in mind that $M^1_v$ is the maximal
class of test functions that works simultaneously for all \modsp s
$\Mmpq $, $1\leq p,q\leq \infty $ and \emph{all $v$-moderate weights
  $m$.}



 Among the modulation spaces occur several  well-known
function spaces : \\
(i) $M^{2}(\Ren)=L^2(\Ren)$. \\
(ii) Weighted $L^2$-spaces: if $\mu _s(x,\o)=\la x\ra^s$, then
$$ M^2_{\mu _s}(\Ren)=L^2_s(\Ren)=\{f\,:\,f(x)\la x\ra^s\in L^2(\Ren)\}.$$
(iii) Sobolev spaces: if  $\ts(x,\o)=\la \o\ra^s$, then
$$ M^2_{\ts}(\Ren)=H^s(\Ren)=\{f\,:\,\hat{f}(\o)\la \o\ra^s\in
L^2(\Ren)\}.$$
 (iv)  Shubin-Sobolev spaces
\cite{Shubin91,BCG02}: If $m_s(z) = (1+|z|)^s = \langle z \rangle ^s $,
then
\begin{eqnarray*}
  M^2 _{m_s}(\rd ) = L^2_s(\Ren)\cap H^s(\Ren)= \Qs(\Ren)
 \end{eqnarray*}
(v) The Schwartz class is related to \modsp s as follows~\cite{GZ01}:
$\cS  (\Ren ) = \bigcap _{s \geq
  0} M^\infty _{\langle \cdot \rangle ^s} (\Ren )\, .$ \\
(vi) The space of tempered distributions~\cite{GZ01}:
\begin{equation}\label{distr}
 \cS ' (\Ren ) = \bigcup _{s \leq 0} M^\infty
_{\langle \cdot \rangle ^s} (\Ren )\,.\end{equation}

\section{Time-frequency tools}

The following    properties of the STFT,  defined in
\eqref{eqi2}, will be used in the sequel. For proofs, see \cite{folland89},
\cite[Ch.~3]{book}.

\begin{lemma}\label{propstft}
Let $f,g,h\in L^2(\Ren)$.  Then we have

(i) Covariance property: if $x,u, \omega , \eta \in \rd $, then 
\begin{equation}
  \label{eq:swi}
  V_g (M_\eta T_u f ) (x,\omega ) = e^{-2\pi i u (\omega -\eta )  } V_gf ( x-u, \omega
  -\eta) \, .
\end{equation}
and thus $|V_g(\pi (y))(z)| = |V_gf(z-y)|$ for $y,z \in \rdd $. 

(ii) Inversion formula:
\begin{equation}
\label{inv} \,\int_{\Renn}  V_gf (x,\om )M_\om T_xh\,dxd\om =\la
h,g\ra f \, ,
\end{equation}
where the integral is to be understood in the weak sense.

\end{lemma}

Fix  a  non-zero window $\gamma$, we define an  operator
 $V_\gamma ^* $ mapping   a function  $F$   on $\rdd$ to a function or distribution on $\rd$  by
\begin{equation}\label{adj}
V^\ast_\gamma F=\int_{\rdd}F\phas M_\o T_x\gamma\,dx\,d\o.
\end{equation}
If  $\gamma , f \in \lrd $ and $F\in L^2(\rdd )$, then $V_\gamma
^*$ is simply the adjoint of $V_\gamma$ in the sense that
$$\la V^\ast_\gamma F, f\ra= \la F, V_\gamma f\ra .$$
 The following
proposition~\cite[Prop.~11.3.2]{book} concerns the boundedness of
the operator $V^\ast_\gamma$ and  will be crucial for the analysis of
localization operators. 

\begin{proposition}\label{Vf*} Let $m\in \cM _v$, $\gamma \in \mvv $,
  $1\leq p,q \leq \infty $. Then   $V_\gamma ^*$ 
maps $L^{p,q} _m (\rdd )$ onto $\Mmpq $
and satisfies
\begin{equation}\label{V*}
\|V_\gamma ^* F\|_{\Mmpq}\lesssim \|V_{\vf _0}
\gamma\|_{L^1_v}\|F\|_{L^{p,q}_m} = \|\gamma \|_{M^1_v}
\|F\|_{L^{p,q}_m} \, .
\end{equation}
 In particular, if $F=V_g f$,
then the inversion formula holds in $\Mmpq$:
\begin{equation}\label{V*2}
f=\frac1{\la\gamma,g\ra}\int_{\rdd}V_g f\phas M_\o T_x \gamma\,dx
d\omega \, .
\end{equation}
\end{proposition}

In this notation, we can write a localization operator informally as
\begin{equation}\label{locaj}
\gaw f =V^*_{\f_2}(a \,V_{\f_1 }f).
\end{equation}

\begin{lemma}\label{Gau}
Let $v$ be any submultiplicative weight  and $\win(t)=e^{-\pi
t^2}$. Then, for all multi-indices $\a,\b\geq 0$,
\begin{equation}\label{*Gau}
\|V_{\win} ( X^\b \partial^\a\win)\|_{L^1_v}=\int_\rdd |V_{\win}
(X^\b \partial^\a  \win)\phas|v\phas\,dxd\o<\infty.
\end{equation}
\end{lemma}
\begin{proof}
We use the (multivariate) Hermite functions $h_\gamma $  to show that
the  integral \eqref{*Gau} is 
finite.  Since $X^\beta\partial^\alpha \f_0$ can be written in the
form $p \f _0$ where $p$ is a polynomial of degree $|\alpha | + |\beta
|$, $X^\beta\partial^\alpha \f_0$ is a finite sum of Hermite functions
$X^\beta\partial^\alpha \f_0=\sum_{|\gamma| \leq |\alpha |+ |\beta |}
c_\gamma  h_\gamma$. By a well-known formula (called the
``Laguerre connection'' in \cite{folland89}), the STFT of
$h_\gamma $ is  $V_{{\win}} h_\gamma (x,\o)=L_{0,\gamma}(\pi (x^2+\o^2))(x,\o) e^{-\pi
(x^2+\o^2)/2},$ where $L_{0,\gamma}$ is a  Laguerre polynomial of
degree $|\gamma |$. Consequently, 
$V_{{\win}}(X^\b\partial^\a{\win})=|P(x,\omega )| e^{-\pi (x^2
  +\omega ^2)/2}$ for some polynomial of degree at most $|\alpha | +
|\beta |$. 
 Since a submultiplicative weight $v$ grows at most
exponentially (i.e., $v\phas\leq\cO(e^{c(|x|+|\o|)}))$, the
convergence  of the integral \eqref{*Gau} is guaranteed.
\end{proof}

Lemma~\ref{Gau} helps us  prove the
boundedness of the operators $ X^\beta\partial^\alpha$ from
$M^1_{v(\la\cdot\ra^N\otimes\la\cdot\ra^N)}$ into
$M^1_{v(\la\cdot\ra^{N-|\b|}\otimes\la\cdot\ra^{N-|\a|})}$.
A similar result can be found in \cite{Toftweight}.

\begin{lemma} \label{finestra}
Let  $g\in M^1_{v(\la\cdot\ra^N\otimes\la\cdot\ra^N)}$ and 
$|\a|,|\beta|\leq N$. Then 
\begin{equation}\label{V*5}
\| X^\beta \partial^\alpha
g\|_{M^1_{v(\la\cdot\ra^{N-|\b|}\otimes\la\cdot\ra^{N-|\a|})}}\lesssim
\|g\|_{M^1_{v(\la\cdot\ra^{N}\otimes\la\cdot\ra^{N})}}.
\end{equation}
\end{lemma}
\begin{proof}
We use the following  algebraic formula~\cite[Lemma~11.2.1]{book} to  interchange
the operators $\partial^\alpha X^\beta$ with \tfs s:
\begin{equation}\label{V*4}
\partial^\alpha X^\beta(M_\o T_x h)=\sum_{\delta_1\leq\alpha}\sum_{\delta_2\leq\beta}\binom{\a}{\delta_1}\binom{\beta}
{\delta_2}x^{\delta_2}(2\pi i \o)^{\delta_1}M_\o T_x
(\partial^{\alpha-\delta_1}X^{\beta-\delta_2}h),
\end{equation}
 for all $\phas\in\rdd$ and $\a,\beta$ with $|\a|,|\beta|\leq N$.
Substituting~\eqref{V*4} into the STFT of 
$ X^\beta \partial^\alpha g$, we obtain 
\begin{align*}
|V_{\win}(X^\beta \partial^\alpha g)(x,\omega )|&=|\la X^\beta
\partial^\alpha g, M_\o T_x\win\ra|=|\la g,\partial^\alpha X^\beta( M_\o
T_x\win)\ra|\nonumber\\
&\leq
\sum_{\delta_1\leq\alpha}\sum_{\delta_2\leq\beta}\binom{\a}{\delta_1}\binom{\beta}
{\delta_2}|x^{\delta_2}(2\pi i \o)^{\delta_1}|\,|\la g, M_\o T_x
(\partial^{\alpha-\delta_1}X^{\beta-\delta_2}\win)\ra|\\
&\leq (2\pi )^N\la x\ra^{|\b|} \la \o\ra^{|\a|}\!\!
\sum_{\delta_1\leq\alpha}\!\sum_{\delta_2\leq\beta}\!\!\binom{\a}{\delta_1}\binom{\beta}
{\delta_2} \!|\la g, M_\o T_x
(\partial^{\alpha-\delta_1}X^{\beta-\delta_2}\win)\ra|.
\end{align*}
Taking the  $L^1$-norm with weight
$v\, (\la\cdot\ra^{N-|\b|}\otimes\la\cdot\ra^{N-|\a|})$,
the previous estimate yields
$$
\|V_{\win}(X^\beta
\partial^\alpha
g)\|_{L^1_{v(\la\cdot\ra^{N-|\b|}\otimes\la\cdot\ra^{N-|\a|})}}\!\!\lesssim\!\!
\sum_{\delta_1\leq\alpha}\sum_{\delta_2\leq\beta}\!\!\binom{\a}{\delta_1}\!\binom{\beta}
{\delta_2}
\!\|V_{(\partial^{\alpha-\delta_1}X^{\beta-\delta_2}\win)} g
\|_{L^1_{v(\la\cdot\ra^{N}\otimes\la\cdot\ra^{N})}}.
$$
 Now  we apply  Theoerem~\ref{convV*3}(i) to each of the terms on the 
 right-hand side and obtain that 
$$\|X^\beta \partial ^\alpha   g
\|_{M^1_{v(\la\cdot\ra^{N}\otimes\la\cdot\ra^{N})}}\lesssim
\sum_{\delta_1\leq\alpha}\sum_{\delta_2\leq\beta}
\|g\|_{M^1_{v(\la\cdot\ra^{N}\otimes\la\cdot\ra^{N})}}
\, \|\partial^{\alpha-\delta_1}X^{\beta-\delta_2}\win  
\|_{M^1_{v(\la\cdot\ra^{N}\otimes\la\cdot\ra^{N})}} <\infty \, .
$$
The latter expression is finite, since  by assumption   $g \in
M^1_{v(\la\cdot\ra^{N}\otimes\la\cdot\ra^{N})}$ and $\partial ^\alpha
X^\beta \win \in
M^1_{v(\la\cdot\ra^{N}\otimes\la\cdot\ra^{N})} \,  
$ by Lemma \ref{Gau}.
\end{proof}

\rem\ Similarly, one can see that $X^\beta \partial ^\alpha $ maps 
$M^{p,q}_{m(\la\cdot\ra^N\otimes\la\cdot\ra^N)}$ into \\
$M^{p,q}_{m(\la\cdot\ra^{N-|\b|}\otimes\la\cdot\ra^{N-|\a|})}$.


\section{Mapping Properties of Localization Operators between
  Modulation Spaces}\label{mapc1}

We  investigate how  a localization operator  maps \modsp s into
each  others. For a more detailed analysis  of the boundedness
properties of localization operators we refer to \cite{CG02}.

\begin{theorem}\label{Elena1}
 Let $m \in \cM _{v}$, $\mu \in \cM _{w}$. 

 \ni (i) Assume that  $a\in L^\infty_{1/m}(\rdd)$,  $\f_1\in M^1_{vw}$
 and 
 $\f_2\in M^1_{w}$, then the localization operator $\gaw$ is
 bounded  from $M^{p,q}_{\mu m}$ into $M^{p,q}_{\mu}$ with a norm
 estimate 
\begin{equation}\label{weu1}
\|\gaw f\|_{M^{p,q}_{\mu}} \lesssim 
\| \f _1 \| _{M^1_{vw}} \, \| \f_2\|_{M^1_w} \, \|a\|_{L^\infty _{1/m}} \,
\|f\|_{M^{p,q}_{\mu m}}. 
\end{equation}

\ni (ii) If  $a\in L^\infty_m(\rdd)$,  $\f_1\in M^1_{w}$,  and
 $\f_2\in M^1_{ vw}$, then  $\gaw$ is 
 bounded  from $M^{p,q}_{\mu}$ into $M^{p,q}_{\mu m}$, with
\begin{equation}\label{weu2}
\|\gaw f\|_{M^{p,q}_{\mu m}} \lesssim \|\f _1\|_{M^1_w} \, \|
\f_2\|_{M^1_{vw}}\, \|a\|_{L^\infty _m} \, \|f\|_{M^{p,q}_{\mu }}.
\end{equation}
\end{theorem}

\begin{proof}
We  write the localization operator $\gaw$ as in \eqref{locaj}
 and  use Proposition \ref{Vf*}.

\ni {\it{(i)}} Since $\f_2\in M^1_{w}$,  estimate  \eqref{V*}
implies that $V_{\f_2}^*$ is bounded from $L^{p,q}_\mu$ onto
$M^{p,q}_{\mu}$ and, for $f\in M^{p,q}_{\mu m}$, we have
\begin{align*}
\|\gaw f\|_{M^{p,q}_{\mu}}&=\|V^*_{\f_2}(a \,V_{\f_1
}f)\|_{M^{p,q}_{\mu}}\\
& \lesssim \| \f_2\|_{M^1_w}\| a \,V_{\f_1 }f\|_{L^{p,q}_{\mu}}.
\end{align*}
Since by assumption $|a\phas|\leq \|a\|_{L^\infty_{1/m}}\,m\phas $ for
$\phas\in\rdd,$  the $L^{p,q}_{\mu}$-norm of the product
$a\,V_{\f_1} f$ can be controlled in the following way:
$$\| a \,V_{\f_1
}f\|_{L^{p,q}_{\mu}}\leq \|a\|_{L^\infty_{1/m}}\|m V_{\f_1}
f\|_{L^{p,q}_{\mu}}=\|a\|_{L^\infty_{1/m}}\| V_{\f_1}
f\|_{L^{p,q}_{\mu m}}\leq
\|\f _1 \|_{M^1_{vw}}\, \|a\|_{L^\infty_{1/m}}\, \|f\|_{M^{p,q}_{\mu m}}.$$ In the last
equivalence we have used $\f_1\in M^1_{vw}$ and Theorem~\ref{convV*3} $(ii)$.

\ni{\it{(ii)}} is similar. Since   $\f_2\in M^1_{vw}$, $V_{\f_2}^*$ is bounded from $L^{p,q}_{\mu
m}$ onto $M^{p,q}_{\mu m}$   by Proposition  \ref{Vf*}. Hence, for $f\in M^{p,q}_{\mu }$, we
have
\begin{align*}
\|\gaw f\|_{M^{p,q}_{\mu m}}&=\|V^*_{\f_2}(a \,V_{\f_1
}f)\|_{M^{p,q}_{\mu m}}\\
& \lesssim \| \f_2\|_{M^1_{ vw }}\| a \,V_{\f_1
}f\|_{L^{p,q}_{\mu m}}\\
&\leq \| \f_2\|_{M^1_{vw}}\|a\|_{L^\infty_{m}}\,\|(1/m)V_{\f_1 }
f\|_{L^{p,q}_{\mu m}}\\ 
&=\| \f_2\|_{M^1_{vw }}\, \|a\|_{L^\infty_{m}}\, \|V_{\f_1
} f\|_{L^{p,q}_{\mu}} \lesssim \|\f _1 \|_{M^1_w}\,  \| \f_2\|_{M^1_{w 
  v}} \, 
\|a\|_{L^\infty_{m}}\,\|f\|_{M^{p,q}_{\mu}}.
\end{align*}
\end{proof}

\rem\  Note once again the role of the
windows. In (i) the source space $M^{p,q}_{\mu m}$ is measured by the
window $\f _1$, which by Theorem~\ref{convV*3} has to be in
$M^1_{vw}$. For the target space $M^{p,q}_\mu $ we use $\f _2$, thus
the required condition is $\f _2 \in M^1_w$. 

\vspace{3 mm}

The following result will be used in the proof of the Fredholm
property of $\gaw$ (Corollary \ref{fine}). For 
weights  of polynomial growth,  Lemma \ref{compact}  was already
observed in~\cite[Lemma 3.8]{PaoloToft}.

\begin{lemma}\label{compact} Let $m \in \cM _{v}$, $\mu \in \cM
  _{w}$. If $a\in L^\infty(\rdd)$ with compact support, and 
$\f_1\in M^1_{v}$, $\f_2\in M^1_{w}$, then the localization operator
$\gaw$ is compact  from $M^{p,q}_{ m}$ into $M^{p,q}_{\mu}$, for
$1\leq p,q\leq\infty$. 
\end{lemma}

\begin{proof} We show that $\gaw$ is bounded and compact  from $M^\infty_{1/v}$ into $M^1_w$. 
 The compactness of $\gaw$ from $M^{p,q}_{ m}$ into $M^{p,q}_{\mu}$
 then follows from the  continuous  embeddings    $M^{p,q}_{m}
 \hookrightarrow    M^\infty_{1/v}$ and  $M^1_w   \hookrightarrow M^{p,q}_{\mu}$. 
We denote the compact support of $a$ by $K\subseteq \rdd $. 

 Let $f\in M^\infty_{1/v}$, then
\begin{eqnarray*}
\|\gaw f\|_{M^{1}_{w}} &=& \| V^*_{\f_2}(a V_{\f_1}f) \|_{M^{1}_{w}}\lesssim \|a V_{\f_1}f\|_{L^{1}_{w}}\\
&=&\int_{K}|a(z)|\,|V_{\f_1}f(z)|w(z)\,dz\\
&\lesssim& \sup_{z\in K} |V_{\f_1}f(z)|\frac{1}{v(z)}\,
\int_{K}|a(z)|v(z)w(z)\,dz\\ 
&\lesssim& \|f\|_{M^\infty_{1/v}} \, , 
\end{eqnarray*}
and so $\gaw $ is bounded. 

Let  $f_n  \in M^\infty _{1/v}$ be a bounded sequence that converges weak$^*$ to some $f
\in M^\infty _{1/v} $. This is equivalent to saying that $V_{\f _1}
f_n $ 
converges uniformly on   compact sets of $\rdd $ to $V_{\f _1 } f$,
e.g., by~\cite[Thm.~4.1v]{fg89jfa}. Consequently,
$$
\| \gaw (f_n -f) \|_{M^1_w} \lesssim \sup_{z\in K} |V_{\f_1}(f_n-f)(z)|\frac{1}{v(z)}\,
\int_{K}|a(z)|v(z)w(z)\,dz \to 0 \,  ,
$$
and this property implies that $\gaw $ is compact from $M^\infty
_{1/v}$ to $M^1_w$. 
\end{proof}

For the analysis of the remainder, we will use the  well-known fact
that every linear  operator $A:\,\cS\rightarrow 
\cS'$ can be written in the form 
$Af=V^*_{\f_2}(T \,V_{\f_1}f)$ for a suitable  integral
operator $T$, see e.g. \cite{folland89,g-heil99}. The boundedness
properties of  operators given in this form are  derived as  in
Theorem~\ref{Elena1}. 
\begin{lemma}\label{TTTT}
Let $\f_1,\f_2 \in M^1_w(\rd)$ and $T$ be the integral operator with kernel $K$ on $\R^{4d}$ (acting on functions $F$ on $\rdd$) defined by $TF(y)=\intrdd K(y,z)F(z)\,dz$   and define 
the operator $A$ (acting on functions $f$ on $\rd$) by 
\begin{equation*}
Af=V^*_{\f_2}(T \,V_{\f_1}f).
\end{equation*}
If $\mu\in\cM_w$ and $T$ is bounded on $L^{p,q}_{\mu}(\rdd)$, $1\leq p,q\leq \infty$, then
$A$ is bounded on $M^{p,q}_{\mu}(\rd)$.
\end{lemma}
\begin{proof}
Since $T$ is bounded by assumption and $V^*_{\f_2}$ is bounded by Proposition \ref{Vf*}, we obtain that 
\begin{align*}
\|A f\|_{M^{p,q}_{\mu}}&=\|V^*_{\f_2}(T \,V_{\f_1
}f)\|_{M^{p,q}_{\mu}}
 \leq \| \f_2\|_{M^1_w}\| T(V_{\f_1 }f)\|_{L^{p,q}_{\mu}}\\
& \leq C_T\| \f_2\|_{M^1_w}\| V_{\f_1 }f\|_{L^{p,q}_{\mu}}
 \lesssim C_T\| \f _1 \|_{M^1_w} \, \| \f_2\|_{M^1_w}\, \|
 f\|_{M^{p,q}_{\mu}}, 
\end{align*}
where the constant $C_T$ is the operator norm of $T$ on $L^{p,q}_{\mu}(\rdd)$.
\end{proof}


\section{Treatment of  Remainder Term}\label{c2}

 The mapping properties of localization operators, as  studied in the
 previous section, enable us to understand the left-hand side in the
 expansion formula ~\eqref{c1}. We now turn to
 the investigation of the remainder term $E_N$.

An explicit formula for the remainder and its Weyl  symbol was
derived in \cite{CR03}. Here we give a different treatment that
leads to estimates for the operator norm of $E_N$. 

We recall   the form  of the remainder
term from \cite{CR03} and the introduction~\eqref{c1}. Let $b_N (y,z) $ be  given by
\begin{equation}
  \label{eq:c11}
b_N(y,z)=N\sum_{|\a|=N}\int_0^1(1-t)^{N-1}\partial^\a
b(y+t(z-y))\,dt\frac{(z-y)^\a}{\a!} \, ,  
\end{equation}
then   the remainder in \eqref{c1} is given by the formula
\begin{equation}
  \label{eq:neu10}
  E_N f = \intrdd \intrdd a(y) b_N (y,z) \, V_{\f _3} f (z) \langle
  \pi (z) \f _4, \pi (y) \f _1\rangle \pi (y) \f _2 \, dz dy \, .
\end{equation}

If we introduce the integral operator $T$ with kernel
\begin{equation}
  \label{eq:68}
  K(y,z) = a(y) b_N (y,z) \,  \langle
  \pi (z) \f _4, \pi (y) \f _1\rangle  \, ,
\end{equation}
i.e., $T H (y)= \intrdd K(y,z) H(z) \, dz \, ,$
 then $E_N $ can be written formally as
\begin{equation}
  \label{eq:670}
  E_N f = V_{\f_2}^*\, (T\,V_{\f_3 }f).
\end{equation}

After these preparations we can formulate the following result for the boundedness and
compactness of the remainder term.

\begin{theorem}\label{resto}
 Let $m\in \cM _v $, $\mu \in \cM _w$,  $1\leq p,q\leq
 \infty$ and assume that $\f _1,\f_4 \in
 M^1_{ w v\langle \cdot \rangle ^N}(\rd )$, and  $\f _2, \f_3 \in
 M^1_{ w }(\rd )$.
 
 \ni
 (i) If $a \in L^\infty _{1/m}(\rdd )$ and $\partial ^\alpha b \in L^\infty
_m (\rdd )$ for $|\alpha| =N$, then $E_N$ is bounded on $M^{p,q}
_\mu $ with the following  estimate
\begin{equation}
  \label{eq:neu2}
  \|E_N f\| _{M^{p,q}_\mu} \lesssim \!\|a\|_{L^\infty _{1/m}} \! \Big(\!\sum
  _{|\alpha| =N}\! \frac{1}{\alpha !} \|\partial ^\alpha b \|_{L^\infty
    _m} \!\Big)\! \|\f _1\|_{M^1_{vw\langle \cdot \rangle ^N}} \!  \| \f
  _2 \|_{M^1_w}\!  \|\f _3 \|_{M^1_w} \!\|\f
  _4\|_{M^1_{vw\langle \cdot \rangle ^N}}\!\!  \| f\| _{M^{p,q}_\mu}.
\end{equation}
(ii) If  $\partial ^\alpha b \,  m \in L^0 (\rdd )$,
then $E_N$ is compact on  $M^{p,q}_\mu $.
\end{theorem}

\begin{proof} $(i)$ It suffices to show that the integral operator $T$
  defined by the kernel $K$ in
\eqref{eq:68}  is bounded on $L^{p,q}_\mu (\rdd )$. Then $E_N$ 
is bounded on $M^{p,q}_\mu $ by Lemma \ref{TTTT}.

{\bf Step $1$.}   An estimate for the kernel $K(y,z)$.

 \ni By
assumption $|a(y)| \leq \|a\|_{L^\infty _{1/m}} \, m(y) $ and
$|\partial ^\alpha b(y)| \leq \| \partial ^\alpha b \|_{L^\infty
_m} \, m(y)\inv$. Therefore we find that
\begin{eqnarray*}
  |a(y) b_N(y,z) | &=& |a(y)| \, N\sum_{|\a|=N}\left|
  \int_0^1(1-t)^{N-1}\partial^\a b(y+t(z-y))\,dt\frac{(z-y)^\a}{\a!} \right|   \\
&\leq & N \|a\|_{L^\infty _{1/m}}\!\! m(y) \la z-y\ra
^N\sum_{|\a|=N}\frac{1}{\a!}\!
   \|\partial^\a b\|_{L^\infty _m}   \! \int _0 ^1\!\!
    m(y+t(z-y))\inv dt  \\
&\leq &  C N \|a\|_{L^\infty _{1/m}}\!\! m(y) \, \la z-y\ra
^N\!\!\sum_{|\a|=N}\frac{1}{\a!}\!
   \|\partial^\a b\|_{L^\infty _m}   \! \int _0^1 \!\! m(y)\inv
v(t(y-z))   dt   \\
&\lesssim &  \|a\|_{L^\infty _{1/m}} \sum_{|\a|=N}\frac{1}{\a!}
   \|\partial^\a b\|_{L^\infty _m}  v(y-z) \langle y-z \rangle
^N
\end{eqnarray*}
where in the last inequality we have used  property
\eqref{*} of the weight $v$.  Consequently, if we set
$C_{a,b,N}:=\displaystyle{\|a\|_{L^\infty _{1/m}}
\sum_{|\a|=N}\frac{1}{\a!}
   \|\partial^\a b\|_{L^\infty _m}}<\infty$,
the integral kernel $K$ is dominated by a convolution kernel in
the sense that
\begin{equation}
  \label{eq:690}
  |K(y,z) | \lesssim C_{a,b,N} v(y-z) \langle y-z \rangle ^N |V_{\f _1} \f _4
  (y-z)| \, .
\end{equation}

{\bf Step $2$.} Boundedness of the integral operator $T$ on
$L^{p,q}_\mu (\rdd )$. 

\ni
By \eqref{eq:690} we obtain that
\begin{eqnarray*}
|T H (y)| &=& \left|\intrdd K(y,z) H(z) \, dz \right| \\
& \leq & C_{a,b,N}    \intrdd v(y-z) \langle y-z
\rangle ^N |V_{\f _1} \f _4   (y-z)| \,|H(z) | \, dz \, \\
&=& C_{a,b,N} \big(|H| \ast (v \langle \cdot \rangle ^N |V_{\f _1} \f _4
|)\big) (y) \, ,
  \end{eqnarray*}
If $H\in L^{p,q}_\mu (\rdd )$, then  Young's inequality with weights
(e.g.~\cite[Prop. 11.1.3.]{book}) implies that 

\begin{eqnarray*}
  \|T H \|_{L^{p,q}_\mu } &\leq & C_{a,b,N} \||H| \ast (v \langle \cdot \rangle
  ^N |V_{\f _1} \f _4 |) \|_{L^{p,q}_\mu } \\
  &\lesssim  &  \|H\|_{L^{p,q}_\mu   }\, \| |V_{\f _1} \f _4 | \, v \langle
  \cdot \rangle   ^N  \|_{L^1_w} \\
&=  &  \|H\|_{L^{p,q}_\mu   }\, \| V_{\f _1} \f _4 \|_{L^1_{vw
\langle \cdot \rangle   ^N}}.
\end{eqnarray*}

{\bf Step $3$.} The boundedness of the remainder $E_N $ is now a
consequence of Lemma \ref{TTTT}.
Accordingly, we need $\f _2 , \f _3 \in M^1_w$, whereas by Theorem~\ref{convV*3} $ \| V_{\f _1}
\f _4 \|_{L^1_{vw  \langle
  \cdot \rangle   ^N}} \lesssim \|\f _1 \|_{M^1_{vw \langle
  \cdot \rangle   ^N}} \, \|\f _4 \|_{M^1_{vw \langle  \cdot \rangle   ^N}} $.  By keeping track of
all the  constants, we find  estimate~\eqref{eq:neu2}. 

\vspace{3 mm}

 \ni (ii) \emph{Compactness of $E_N$.}

{\bf Step $4$}. We first show that if  $b$ has compact support
then $E_N$ is a compact operator    from $M^\infty_\mu$ into
$M^1_{\mu}$. So assume that  $\supp b \subseteq K$ for some
compact set $K\subseteq \rdd $. Similar arguments as in item $(i)$
give the boundedness of the integral operator $T$ from
$M^\infty_\mu$ to $M^1_\mu$. Since by
assumption $\partial^\a b\in L^\infty_m$, we have
$$ |\partial^\a b(z)|\leq\|\partial^\a
b\|_{L^\infty_m}m(z)^{-1}\chi_K(z).$$ 
As in Step $1$ above we derive the kernel estimate
\begin{eqnarray*}
|K(y,z)|& \lesssim & \|a\|_{L^\infty_{1/m}}\sum_{|\a|=N}\frac1{\a!}\|\partial^\a
b\|_{L^\infty_m} \int_0^1  \, v(t(z-y)) \chi_K(y+t(z-y))dt\cdot\\
& \, &\qquad\qquad\qquad \cdot\,\la z-y\ra^{N}
|V_{\f_4}\f_1(z-y)|\,.
\end{eqnarray*}

Let $H\in L^\infty_\mu(\rdd)$. To estimate $TH$ we apply the
preceding computation  and also use the $w$-moderateness of $\mu
$. This yields the majorization
\begin{eqnarray}\|TH\|_{L^1_{\mu}}&\leq&
\intrdd\left(\intrdd |K(y,z)| |H(z)|\,dz\right)\mu(y)\,dy\nonumber\\
&\lesssim& \intrdd\intrdd\int_0^1 v(t(z-y))\chi_K(y+t(z-y))| \, \,
\la
z-y\ra^{N}H(z)\cdot\label{stima}\\
&\,&\qquad\qquad \cdot\,|V_{\f_4}\f_1(z-y)| \, \mu(z)w(z-y)\, dz
dy dt \,.\nonumber
\end{eqnarray}
Performing the change of variables $z'=z-y$, $y'=y+t(z-y)$ so that  $dydz =
dy' dz'$, and using \eqref{*}, we get
\begin{eqnarray}\label{es1}\|TH\|_{L^1_{\mu}}&\lesssim&\|H\|_{L^\infty_\mu}\intrdd\intrdd\int_0^1
v(t z')\chi_K(y')\la z'\ra^N |V_{\f_4}\f_1(z')|w(z')\, dz' 
dy' dt\,\nonumber\\
&\lesssim&\|H\|_{L^\infty_\mu} \intrdd \chi_K(y')\,dy' \intrdd v(
z')w(z')\la z'\ra^N |V_{\f_4}\f_1(z')|\,dz'\nonumber\\
&=&\|H\|_{L^\infty_\mu}\,
meas(K)\,\|V_{\f_4}\f_1\|_{L^1_{vw\la\cdot\ra^N}}
\end{eqnarray}

To prove the compactness, we need to show that every bounded
sequence $f_n\in M^\infty_\mu$ possesses  a subsequence $f_{n_k}$
such that $E_N f_{n_k}$ converges in the $M^1_{\mu}$-norm. So, assume that
$f_n\in M^\infty_\mu = (M^1_{1/\mu })^*$ and $\|f_n \|_{M^\infty
_\mu } \leq 1 $ for all $ n\in \bN $. By the Theorem of
Alaoglu-Bourbaki there is a subsequence $f_{n_k}$
that converges in the $w^*$-topology to some $f\in M^\infty_\mu$.
After replacing $f_{n_k}$ by $f_{n_k} -f$, we may
 assume without loss of generality that  $f_{n_k} \stackrel{w^*}{\to} 0$.
Since $\f_3\in M^1_{w}(\rd)\subseteq M^1_{1/\mu}(\rd)$,  
the $w^*$-convergence  implies in particular that
$$(V_{\f_3} f_{n_k})(z)=\la f_{n_k},\pi(z)\f_3\ra \rightarrow
0, \quad \quad \forall z\in \rdd \, .$$ 
The estimate \eqref{stima} yields  the majorization 

\begin{eqnarray*}\|E_N f_{n
_k}\|_{M^1_\mu}&\lesssim& \intrdd\intrdd\int_0^1 v(t(z-y))\chi_K(y+t(z-y))| \, \,
\la
z-y\ra^{N}V_{\f_3}f_{n_k}(z)\cdot\\
&\,&\qquad\qquad \cdot\,|V_{\f_4}\f_1(z-y)| \, \mu(z)w(z-y)\, dz
dy dt \,.\nonumber
\end{eqnarray*}
Since $V_{\f_3}f_{n_k}(z)\mu(z)\leq \|f_{n_k}\|_{ M^\infty_\mu}\leq 1$ for all $z\in\rdd$, the previous integral is
bounded (uniformly in $n_k$)  by the  function 
 $$v(t(z-y)) \, \chi_K(y+t(z-y)) \,  \la
z-y\ra^{N} \, |V_{\f_4}\f_1(z-y)| \, w(z-y).$$
 This function  is  
integrable  on $\rdd \times \rdd \times [0,1]$ by~\eqref{es1}.
Thus the hypotheses of Lebesgue's Theorem on 
 dominated convergence   are satisfied, and we
conclude that 
 $\|E_N f_{n_k}\|_{M^1_\mu } \rightarrow 0$,  as desired. So $E_N$ is
 compact  from $M^\infty_\mu$ into $M^1_{\mu}$.

\ni {\bf Step $5$}. Since the embeddings $\mpq _\mu
\hookrightarrow M^\infty_\mu$ and
$M^1_{\mu}\hookrightarrow\mpq_{\mu}$ are continuous, it follows
that $E_N$ is compact from $\mpq_{\mu}$ into $\mpq_{\mu}$.

\ni {\bf Step $6$}. If $(\partial^\a b) m  \in L^0$ for all
$|\a|=N$, there exist  sequences of functions $\rho^\a_n$ with
compact support such that
$$\|(\partial^\a b) m-(\rho^\a_n)  m\|_{L^\infty}\rightarrow 0,\qquad \forall \a,\,|\a|=N.$$

Let $E_N^n$ denote  the  operator obtained by replacing the
derivatives  $\partial^\a
b$ by their approximations  $\rho^\a_n$ in~\eqref{eq:neu10}. Then by
Step 
$5$  $E_N^n: \mpq_{\mu}\rightarrow \mpq_{\mu}$ is a compact
operator,   and  we have
$$\|E_N-E^n_N\|_{\mpq_{\mu}\rightarrow\mpq_{\mu}}\lesssim \sum_{|\a|=N}\frac{\|\partial^\a b- \rho^\a_n
\|_{L^\infty_m}}{\a!}\rightarrow 0\, .$$
  Being the limit of compact
operators, $E_N$  is also  compact  from $\mpq_{\mu}$
into $\mpq_{\mu}$.
\end{proof}

\section{The Main Theorem -- Symbolic Calculus for Localization
  Operators}

We now have all the pieces in place to prove the validity of the
composition formula for localization operators
\begin{equation}\label{c1a1}
\gaw A^{\f_3,\f_4}_b= \sum_{|\a|=0}^{N-1} \frac{(-1)^{|\a|}}{\a!}
A^{\Phi_\a,\f_2}_{a\partial^{\a} b}+E_N
\end{equation}
  on  modulation spaces. According to \cite{CR03} the  windows $\Phi_\a$ are
   given explicitly  by the   formula
\begin{equation}\label{newwind}
\Phi_\a=\frac1{(2\pi i)^{|\a_2|}}\sum_{\b\leq\a}
\binom{\a}{\b}(-1)^{|\b_1|}\la \f_3,
X^{\a_1-\b_1}\partial^{\a_2-\b_2}\f_4\ra
X^{\b_1}\partial^{\b_2}\f_1 \, .
\end{equation}

In order to treat \eqref{c1a1} on \modsp s, we need to verify that
(a) all occurring windows are in the correct window classes,  
(b) all the terms are well-defined bounded operators on \modsp s, 
and (c) the equality \eqref{c1a1} holds true. The proof consists in
the careful combination of the main results of
Sections~\ref{mapc1} and~\ref{c2}.

The generality of the hypotheses on the admissible  weight functions,
windows, and symbols makes the conditions somewhat technical and
cumbersome. At a first
reading one may assume that all weights are of at most  polynomial
growth. Then it suffices to take all windows $\vf _j , j=1, \dots , 4$,
in $\cS (\rd )$ (whereas the more general conditions in the following
theorem permit to choose them in much larger \modsp s).


\begin{theorem}\label{mainA}
Consider  the following set of hypotheses:

(i) Weights: $m\in \cM _v,  \mu\in\cM _w$.

(ii) Windows:
  $\f_1,\f_4 \in M^1_{vw(\wei^N\otimes\wei^N)}$,  $\f_2,\f_3  \in
M^1_{w}$.

(iii) Symbols:
 $a\in L^{\infty}_{1/m}(\rdd)$,
$\partial^\a b \in L^{\infty}_{m}(\rdd)$ for all $|\a|\leq N$. \\

If hypotheses  (i) --- (iii) hold, then the symbolic calculus for localization
operators is valid in the following sense:

\ni
 (i) For $|\alpha| < N$  the  window $\Phi_\a$ on
$\rd$ defined by~\eqref{newwind}
belongs to $M^1_{vw}$.\\
\ni (ii) The composition formula 
given in \eqref{c1a1}
 is well-defined on
every $\mpq_\mu$, for $1\leq p,q\leq\infty$. This means that  the
product on the left-hand side and  all operators   on the
right-hand side are bounded  from  $\mpq_{\mu}$ to $\mpq_{\mu} $,
and
 \eqref{c1a1}  is valid as  an identity of operators.

(iii) If, in addition,
\begin{equation}
  \label{decay}
   (\partial ^\alpha b )\, m \in L^0 (\rdd ) \, ,\quad \forall \,\alpha,
   |\a|=N,
\end{equation}
    then the  remainder
 $E_N$   is compact on $\mpq_{\mu} $.
\end{theorem}

\begin{proof}
\emph{Boundedness of  the operators on the left-hand side.}
We apply  Theorem~\ref{Elena1} whose assumptions are tailored for this
purpose. By \eqref{weu2} $A_b^{\f_3,\f_4}$ is 
 bounded from $\mpq_\mu$ into $\mpq_{ \mu m}$, and by \eqref{weu1}
 $\gaw$ 
 is  bounded from $\mpq_{ \mu m}$ into $\mpq_\mu$. 
Thus  product   $\gaw  A_b^{\f_3,\f_4}$ is therefore bounded 
 on $\mpq_{ \mu }$ with an operator norm not exceeding $
\|a\|_{L^\infty _{1/m}} \, \| b \|_{L^\infty _m} \,
  \|\f _1  \|_{M^1_{vw}} \|\f _2  \|_{M^1_{w}} \|\f _3  \|_{M^1_{w}}
  \|\f  _4  \|_{M^1_{vw}}  $.

\vskip0.3truecm\ni
\emph{Boundedness of  the operators on the right-hand side.}
The boundedness of the remainder $E_N$ on $M^{p,q}_\mu $  was stated
and proved in Theorem~\ref{resto}(i). With the additional decay
condition~ \eqref{decay} the compactness of $E_N$ follows from
Theorem~\ref{resto}(ii).  Note that $\langle (x,\omega ) \rangle ^N
\leq \langle x \rangle ^N \, \langle \omega \rangle ^N$ and thus
$M^1_{vw(\wei^N\otimes\wei^N)} \subseteq M^1 _{vw \wei ^N}$. The
slightly stronger condition on $\f _1, \f _4$ is required in the next
step. 


Next we turn to the  sum of localization operators
$A^{\Phi_\a,\f_2}_{a\partial^\a b}$. We first check that each  window $\Phi_\a$
defined in \eqref{newwind} belongs to $M^1_{v w }$. For
$\b = (\b_1,\b_2) $ and $|\b | \leq N$,    Lemma~\ref{finestra} gives
$$
X^{\b_1}\partial^{\b_2}\f_1\in
M^1_{vw(\wei^{N-|\b_1|}\otimes \wei^{N-|\b_2|})}\subseteq
M^1_{vw} \, $$
and likewise $ X^{\a_1-\b_1}\partial^{\a_2-\b_2}\f_4 \in M^1_{vw}
\subseteq L^2$ when $|\beta | \leq |\alpha | \leq N$.
Consequently the
brackets $\la \f_3, t^{\a_1-\b_1}\partial^{\a_2-\b_2}\f_4\ra$ are
  well-defined  as  an inner product in $L^2$, and  $\Phi_\a$, as a
finite  linear combination of $X^{\b_1}\partial^{\b_2}\f_1$,   is then
also  in $M^1_{vw}\subseteq  M^1_w$. 

For  boundedness of $A^{\Phi_\a,\f_2}_{a\partial^\a b}$ on
$M^{p,q}_\mu$ we invoke once more  Theorem~\ref{Elena1}. We know that
$a\partial^\a b \in L^\infty _{1/m} \cdot L^\infty _m  = L^\infty
(\rdd )$  and  $\Phi_\a,\f_2  \in  M^1_w$, hence 
$A^{\Phi_\a,\f_2}_{a\partial^\a b}$ is bounded on $\mpq_\mu$. \vskip0.3truecm

\ni \emph{Formula \eqref{c1a1} is an identity of operators.} 
So far we have shown that all terms in the symbolic
calculus~\eqref{c1a1} are 
well-defined and bounded on $M^{p,q}_\mu $. It remains to show that 
\eqref{c1a1} is actually an identity of operators, under the general
hypothesis stated. 
Regarding this question, we already know from the main result
in~\cite{CR03}
that \eqref{c1a1} holds for windows $\f _j$ in the Schwartz class and
symbols in certain Shubin classes, in particular for symbols  $a,b\in
\cS (\rdd )$. We need to extend the validity of \eqref{c1a1} to
windows in the (possibly larger) modulation spaces and to
\emph{non-smooth} symbols in weighted $L^\infty $-spaces. 

To accomplish this extension, we view each of the terms in the
symbolic calculus as a multilinear form mapping $(a,b,\f _1, \dots ,
\f _4)= (a,b, \vec{\f})$ to one of the operators 
\begin{equation}
  \label{eq:cn1}
  T(a,b,\vec{\f }) = \gaw A_b ^{\f _3, \f _4}, \,\, A^{\Phi _\alpha, \f _2}
  _{a\partial ^\alpha b}, \, \, \mathrm{or}\, \, E_N (a,b, \f) \, .
\end{equation}

  The operator norm of each of  these operators on
$M^{p,q}_\mu $ obeys  an estimate of the form
\begin{eqnarray}
\label{eq:cn2}
\lefteqn{  \|T(a,b,\vec{\f } )\|_{M^{p,q}_\mu  \to M^{p,q}_\mu}
  \lesssim } \\
& \lesssim &  \|a\|_{L^\infty _{1/m}} \, \|\partial ^\alpha b \|_{L^\infty _m} \,
  \|\f _1  \|_{M^1_{vw(\wei^N\otimes\wei^N) }} \|\f _2  \|_{M^1_{w}}
  \|\f _3  \|_{M^1_{w}} \|\f  _4  \|_{M^1_{vw(\wei^N\otimes\wei^N)}}
  \, ,  \notag   
\end{eqnarray}
for suitable $\alpha , |\alpha |\leq N$ as proved in
Theorem~\ref{Elena1} and~\ref{resto}, and~\eqref{eq:neu2}.  

The extension of the symbolic calculus from windows in $\cS (\rd )$ to
windows in the  \modsp s can be done by a routine
density argument, because $ \cS _{\cC } $ is dense in $M^1_v$ for
any submultiplicative weight $v$. Since by \eqref{eq:cn2} each
operator $T(a,b,\vec{\f } )$ is jointly continuous in $\f _j, j=1 , \dots ,
4$, we may choose four sequences $\f _{j,n} \in \cS _{\cC }$
such that $\f _{j,n} \to \f _j$ in the correctly weighted
$M^1$-norm. Then, as $n\to \infty $,  each $T(a,b, \vec{\f _n})$
converges to $ T(a,b, 
\vec{\f}) $ in operator norm. As a consequence,
the symbolic calculus \eqref{c1a1} holds under the general hypotheses
on the windows as stated. 

The extension of the symbolic calculus to non-smooth symbols is more
subtle, because $\cS $ is not norm-dense in $L^\infty _{1/m}$ or
$L^\infty 
_m$. We have to take recourse to a weak-$^*$ approximation argument. 

Given $a \in L^\infty _{1/m}(\rdd )$, we choose a sequence $a_k \in
\cS (\rdd )$ converging weak-$^*$ to $a$, i.e., $\langle a_k , F
\rangle \to \langle a, F \rangle$ for all $F \in L^1_m$. Likewise,
given $b\in L^\infty _m$ with $\partial ^\alpha b \in L^\infty _m$ for
all $\alpha, |\alpha | \leq N$, we may choose a sequence $b_n \in \cS
$, such that $\langle \partial ^\alpha b_n , F \rangle \to \langle
\partial ^\alpha b , F \rangle $ for all $F\in L^1_{1/m}$. This is
always possible by a regularization of $b$, see e.g.~\cite{hormander}.

Next we  show that,  for each operator $T(a,b, \vec{\f})$
in~\eqref{eq:cn1} and all $f\in M^{p,q}_\mu $ and $g\in
M^{p',q'}_{1/\mu }$, we have 
\begin{equation}
  \label{eq:cn3}
  \lim _{k, n \to \infty } \langle T(a_k, b_n, \vec{\f}) f, g
  \rangle  = \langle  T(a, b, \vec{\f}) f, g
  \rangle \, .
\end{equation}

To be specific, we carry out the argument for the remainder
$E_N(a,b,\vec{\f})$. 

Since by Theorem~\ref{resto} we have 
$$
|\langle E_N(a,b)f,g\rangle |\lesssim \|a\|_{L^\infty _{1/m}} \, ,
$$
there exists an $F\in L^1_m$ (depending on all other parameters $b,\f
_j, f, g$), such that 
$$\langle E_N(a,b)f,g\rangle = \langle a, F\rangle \, .
$$
In fact, according to~\eqref{eq:neu10} $F$ is given explicitly by 
$$
F(y) =  \intrdd  b_N (y,z) \, V_{\f _3} f (z) \langle
  \pi (z) \f _4, \pi (y) \f _1\rangle \overline{V_{\f _2} g(y)}\,  dz \, ,
$$
and can be shown directly to be in $L^1_m (\rdd )$  as in Steps~1
and~2 of the proof of Theorem~\ref{resto}.

Consequently, if $a_k \stackrel{w^*}{\to} a$, then 
$$\langle E_N (a_k,b)f,g\rangle = \langle a_k, F\rangle \to \langle
a,F\rangle = \langle E _N(a,b)f,g\rangle  \, .
$$
For the convergence in $b$, let $E_N^\alpha $ be the term in the
remainder  that corresponds to the $\alpha$-th derivative of
$b$. Since as part of the proof of Theorem~\ref{resto} we have shown
that 
$$
|\langle E_N ^\alpha (a,b)f,g\rangle | \lesssim \|\partial ^\alpha b
\|_{L^\infty _m} \, ,
$$
there exists  $G\in L^1 _{1/m}$ (depending on $a,f,g , \f _j$)
such that 
$$
\langle E_N ^\alpha (a,b)f,g\rangle = \langle \partial ^\alpha b , G
\rangle \, .
$$
Consequently, if $\partial ^\alpha b_n \stackrel{w^*}{\to} \partial
^\alpha b$, then 
$$\langle E_N^\alpha (a,b_n)f,g\rangle = \langle \partial ^\alpha b_n,
G\rangle \to \langle \partial ^\alpha b, G\rangle = \langle
E_N^\alpha (a,b)f,g\rangle  \, .
$$
By summing over all $\alpha , |\alpha | =N$,  we have shown
\eqref{eq:cn3} for $E_N$. The weak-$^*$ convergence of 
$\gaw A^{\f _3, \f _4}_b $ and $A^{\Phi _\alpha, \f _2}
  _{a\partial ^\alpha b}$ is shown by exactly the same argument. 

Finally, the weak-$^*$ convergence of each term in the symbolic
calculus \eqref{c1a1}  implies that
it remains valid under weak-$^*$ limits of the 
symbols $a$ and $b$. Consequently, \eqref{c1a1} holds for arbitrary
symbols $a\in L^\infty _{1/m}$ and $b$ with $\partial ^\alpha b \in
L^\infty _m $ for $|\alpha | \leq N$. This completes the proof of the
expansion formula on $M^{p,q}_\mu $. 
\end{proof}

Theorem  \ref{mainA} possesses  a symmetric version that is obtained
by using  the Taylor expansion  for the symbol $a$ instead of
$b$. Precisely,  the following result holds.

\begin{proposition}\label{mainadj}
Under the assumptions of Theorem \ref{mainA}, with the rules of
the symbols $a$ and $b$ interchanged, that is, $b\in L^\infty
_{1/m}$, $\partial^\a a\in L^\infty_m$, $\forall \a\leq N$, we obtain 
the \emph{symmetric} composition formula
\begin{equation}\label{?*e}
\gaw A^{\f_3,\f_4}_b= \sum_{|\a|=0}^{N-1} \frac{(-1)^{|\a|}}{\a!}
A^{\f_3, \Psi _\alpha}_{(\partial^{\a}a) b}+{\widetilde E}_N\, .
\end{equation}
Here 
$$
\Psi _\alpha = \frac{(-1)^\alpha}{(2\pi i)^{|\a_2|}}\sum_{\b\leq\a}
\binom{\a}{\b}(-1)^{|\b_1|}\la \f_2,
X^{\a_1-\b_1}\partial^{\a_2-\b_2}\f_1\ra
X^{\b_1}\partial^{\b_2}\f_4 \, .
$$
and  the remainder ${\widetilde E}_N$ is given by
$${\widetilde E}_Nf=\intrdd \intrdd  a_N(z,y) b(z)V_{\f_3}f(z)\la
\pi(z)\f_4,\pi(y)\f_1\ra\pi(y)\f_2\,dydz \, ,$$ 
with $a_N$ already defined in~\eqref{eq:c11}. 
Formula \eqref{?*e} is a well-defined identity of bounded
operators on $\mpq_\mu$.

(ii) If 
\begin{equation}
  \label{decaya}
   (\partial ^\alpha a )\, m \in L^0 (\rdd ) \, ,\quad \forall \,\alpha,
   |\a|=N,
\end{equation}
    then  $\widetilde{E}_N$   is compact on $\mpq_{\mu} $.
\end{proposition}

\section{Fredholm Property  of Localization Operators}

In the final section we investigate the Fredholm property of
localization operators. As in the theory of partial differential
equations, the construction of a parametrix (a left inverse or a
right inverse modulo regularizing terms) is one of the main
applications of the symbolic calculus.

  Recall that a bounded
 operator $A: B_1 \to B_2$  between two  Banach spaces $B_1, B_2$ is
 called a Fredholm operator if
 $${\rm dim\, Ker\,}A<\infty\quad {\rm\,and}\, \quad{\rm dim\, Coker
}\,A<\infty.$$
Equivalently, $A$ is Fredholm, if there exists a \emph{left parametrix}
$B:B_2 \to B_1$ such that $BA = \mathrm{I}_{B_1} + K_1$ for some compact
operator $K_1: B_1 \to B_1$ and a right  parametrix $C: B_2 \to B_1$ such that $AC =
\mathrm{I}_{B_2} + K_2$ for a compact operator  $K_2: B_2 \to B_2$ \cite{BBR,CR03}.

The expansion formula~\eqref{c1} for $\gaw $ with $N=1$ and $b= 1/a$
yields a 
natural candidate for a parametrix, namely the localization
operator $A^{\f _2,\f _1} _{1/a}$. This idea is formulated
precisely in the next theorem.

\begin{theorem}[{\bf Fredholm property}]\label{F1}
Let $m\in \cM _v, \mu \in \cM _w$. Assume that the symbol  $a$  and the
windows $\f _1, \f _2$ satisfy the following conditions:

\ni
 (i) $|a| \asymp 1/m $, in particular   $a\in L^\infty _m (\rdd )$,    \\

 \ni (ii) $ (\partial_j a)m\in L^0 $  for  $  j=1,\dots,2d$,  and 

\ni (iii) $\f_1,\f_2\in M^1_{w v^2(\wei\otimes\wei)}$ and $|\la
\f_1,\f_2\ra |=1$.

  Then  the  operators
$$\gaw:\,\mpq_\mu\rightarrow\mpq_{\mu m}\qquad \mathrm{ and } \qquad
A^{\f_2,\f_1}_{1/a}:\,\mpq_{\mu m}\rightarrow\mpq_\mu$$ are
Fredholm operators.
\end{theorem}

\begin{proof} \emph{Construction of a  left parametrix for  $\gaw$}.
Once more,  Theorem~\ref{Elena1}(ii) 
 implies   that $\gaw $ is bounded from
  $M^{p,q}_{\mu}$ into $M^{p,q}_{\mu m}$, since $a\in L^\infty
  _{m}(\rdd )$. Likewise, $1/a \in
  L^\infty _{1/m}$ and by  Theorem~\ref{Elena1}$(i)$ 
  $A^{\f_2,\f_1}_{1/a}$ is bounded from 
  $M^{p,q}_{\mu m}$ into $M^{p,q}_{\mu }$.

 Next, we apply  Theorem~\ref{mainA} with  $N=1$ to the 
 product $A^{\f_2,\f_1}_{1/a} \gaw$ and obtain the  composition
 formula 
\begin{equation}\label{Fprod1}
A^{\f_2,\f_1}_{1/a} \gaw=A_1^{\Phi_0,\f_1}+E_1 \, ,
\end{equation}
as  an operator identity  on $\mpq_\mu$. 

According to~\eqref{newwind}, $\Phi _0 $ is given by 
$\Phi _0 = \langle \f _1, \f _2 \rangle \f _2$. Using the 
inversion formula \eqref{inv}, which corresponds to $a \equiv 1$, we 
obtain
\begin{eqnarray*}
A^{\Phi_0,\f_1}_{1}f &=& \int _{\rdd } 1 \, \langle f, \pi (z) \Phi _0
\rangle \f _1 \, dz \\
&=& \langle \f _1, \Phi _0 \rangle f = |\langle \f _1, \f _2 \rangle
|^2 f = f \, .
\end{eqnarray*}

The remainder $E_1$ is compact on $\mpq_\mu$ because of assumption
(ii) on $a$ and assumption (iii) on the windows and
Theorem~\ref{resto}(ii). This means that $A^{\f_2,\f_1}_{1/a}$
is a left parametrix   for  $\gaw$.

\emph{Construction of a right parametrix  of $\gaw$.} The construction of a right
parametrix  is similar. 
Since we apply Proposition~\ref{mainadj} to $M^{p,q}_{\mu m}$ instead of
$M^{p,q}_\mu$, we have to replace the pair of weights $(w, \mu )$ 
in the hypotheses by the pair $(vw, \mu m)$.   
Then, 
  by the symmetric formula~\eqref{?*e}, we obtain the composition formula on $\mpq_{\mu m}$:
\begin{equation*}\gaw
A^{\f_2,\f_1}_{1/a}= A_1^{\Phi _0, \f _2} + \widetilde{E}_1 =
I+\widetilde{E}_1 \, . 
\end{equation*}
Again by Theorem~\ref{resto}(ii)  $\widetilde{E}_1$  is a compact
operator on $\mpq_{\mu m}$, and thus    $ A^{\f_2,\f_1}_{1/a}$ is
the right parametrix for $\gaw$.

Altogether we have shown that 
$\gaw$  is a  Fredholm  operators  between $M^{p,q} _\mu $ and
$M^{p,q}_{\mu m}$. 
Likewise, $  A^{\f_2,\f_1}_{1/a}$ is Fredholm between 
$M^{p,q}_{\mu m}$ and $M^{p,q} _\mu $. 
\end{proof}

Theorem \ref{F1} applies to the standard weight functions, such as 
$(1+|z|^2)^{s/2}$ or 
$(1+|z|)^{s}$, and also to submultiplicative weights satisfying a
condition of the form 
$|\partial_jv(z)|\leq v(z)^{\tau}$ for some $\tau,  0<\tau<1.$ 
However, as written, our main theorem does not seem to work for the
subexponential weights $v(z)=e^{a |z|^b},$ $0<b<1$, because $$
(\partial_j  e^{a |z|^b}) e^{-a |z|^b}=-a |z|^{b-1}\frac {z_j}{|z|}$$ 
possesses a mild singularity at $0$ and is not bounded in a neighbourhood of $0$.
The next result shows that this difficulty  can be circumvented
easily, because the decisive property of $a$ in Theorem~\ref{mainA} is the asymptotic
behavior at $\infty$ and not the local behavior at $0$.

\begin{corollary}\label{fine} The conclusions of Theorem \ref{F1} hold
  if condition $(ii)$ on the symbol $a$ is replaced by the weaker 
  condition:\c 
$(ii')$ there exists a compact set $K\subset\rdd$ such that, 
\begin{equation}\label{comperr}
(1-\chi_K)(\partial_ja) m\in L^0(\rdd),\quad j=1,\dots, 2d.
\end{equation}
\end{corollary}

\begin{proof}
We first find the left-inverse of $\gaw$. 
Let $\psi$ be a test function in $\cD$ such that $\psi(z)=1$ for $z\in K$. Then the product
$A^{\f_2,\f_1}_{1/a} \gaw$ can be recast as
\begin{equation}\label{Ferrprod1}
A^{\f_2,\f_1}_{1/a}
\gaw=A^{\f_2,\f_1}_{1/a}A^{\f_1,\f_2}_{a\psi}+A^{\f_2,\f_1}_{(1/a)\psi}\,
A^{\f_1,\f_2}_{a(1-\psi)}+A^{\f_2,\f_1}_{(1/a)(1-\psi)}A^{\f_1,\f_2}_{a(1-\psi)}
\end{equation}

Since both functions $a\psi$ and $(1/a)\psi$ are bounded with compact
support, Lemma \ref{compact} guarantees that the corresponding
operators are compact on $\mpq_\mu$. Thus, the first two terms of the
right-hand side of \eqref{Ferrprod1} are compact operators 
on $\mpq_\mu$. To treat the third term, we apply the expansion formula
of Theorem \ref{mainA} with $N=1$ and we obtain 
\begin{eqnarray*}
A^{\f_2,\f_1}_{(1/a)(1-\psi)}A^{\f_1,\f_2}_{a(1-\psi)}&=&A_{(1-\psi)^2}^{\Phi_0,\f_1}+E_1\\
&=& I +{\overline{\la \f_1,\f_2\ra}}A_{-2\psi+\psi^2}^{\f_2,\f_1}+E_1.
\end{eqnarray*}
Since $(ii')$ implies that $a(1-\psi)$ satisfies condition \eqref{decaya}, the remainder
$E_1$ is compact on $\mpq_\mu$ whereas  the compactness of the
operator $A_{-2\psi+\psi^2}^{\f_2,\f_1}$ on $\mpq_\mu$ follows again  by
Lemma \ref{compact}. 

We have shown that $A^{\f_2,\f_1}_{1/a} \gaw= I+ R$, for some compact operator $R$ and thus 
have obtained a left parametrix for $\gaw$. The argument for the right parametrix is similar.
\end{proof}

\rem\ It is natural to conjecture that for symbols of the form $a=
1/m$ the localization operator $\gaw $ is an isomorphism between
$M^{p,q}$ and $M^{p,q}_m$. Under very restrictive conditions ($m$
depends only on one variable, or $m$ is hypoelliptic and grows at most
polynomially) this has been proved in
\cite{fei83,PaoloToft} and others. 
Using completely different techniques (Banach algebra methods and
spectral invariance) the following can be shown: assume that $a$ is
measurable and $0< A \leq a(z) \leq B < \infty $ for almost all
$z\in \rdd $ and $\f _1, \f_2 \in M^1_v$, then $\gaw $ is an
isomorphism on $M^{p,q}_m$ for $1\leq p,q\leq \infty $ and all $m\in
\cM _v$~\cite{grochesi}. In contrast to existing results,    no smoothness
conditions are required on the symbol $a$   in this case.

\section*{Acknowledgements}
The authors would like to thank Prof. Luigi Rodino for fruitful
conversations and comments. This work was completed while both authors
were visiting the Erwin Schr\"odinger Institute in Vienna. Its
hospitality and great working conditions are gratefully acknowledged. 


\end{document}